\documentclass[12pt, a4paper, reqno, oneside]{article}

\usepackage{amssymb,amsmath}
\DeclareUnicodeCharacter{0335}{}
\usepackage{tikz-cd}
\usepackage{authblk} 

\usepackage{tikz}
\usepackage{titlesec} 
\usetikzlibrary{babel}
\usepackage{adjustbox}
\usepackage{url}
\usepackage{xcolor}
\usepackage{eso-pic}
\usepackage[T1]{fontenc}
\usepackage[utf8]{inputenc}
\usepackage[english]{babel}
\usepackage{csquotes}
\usepackage{latexsym}
\usepackage{euscript}
\usepackage{mathtools}
\usepackage[all]{xy}
\usepackage{mathrsfs}
\usepackage[percent]{overpic}
\usepackage{setspace}
\usepackage{graphicx}
\usepackage{amsthm}
\usetikzlibrary{calc}
\usetikzlibrary{decorations.pathmorphing}

\tikzset{curve/.style={settings={#1},to path={(\tikztostart)
    .. controls ($(\tikztostart)!\pv{pos}!(\tikztotarget)!\pv{height}!270:(\tikztotarget)$)
    and ($(\tikztostart)!1-\pv{pos}!(\tikztotarget)!\pv{height}!270:(\tikztotarget)$)
    .. (\tikztotarget)\tikztonodes}},
    settings/.code={\tikzset{quiver/.cd,#1}
        \def\pv##1{\pgfkeysvalueof{/tikz/quiver/##1}}},
    quiver/.cd,pos/.initial=0.35,height/.initial=0}

\tikzset{tail reversed/.code={\pgfsetarrowsstart{tikzcd to}}}
\tikzset{2tail/.code={\pgfsetarrowsstart{Implies[reversed]}}}
\tikzset{2tail reversed/.code={\pgfsetarrowsstart{Implies}}}
\tikzset{no body/.style={/tikz/dash pattern=on 0 off 1mm}}

\usepackage[colorlinks=true,linktocpage=true, pagebackref=false, citecolor=black,linkcolor=black]{hyperref}
\normalbaselines
\makeatletter

\def\ps@myfancy{\let\@mkboth\markboth
 \def\@evenhead{\vbox{\hsize\textwidth 
 \hbox to \textwidth{\sf\mdseries\thepage 
 \rule[-.6ex]{0mm}{2mm} \hfill\sf\large\leftmark}
 \vskip 1pt \hrule}}
 \def\@oddhead{\vbox{\hsize\textwidth 
 \hbox to \textwidth{{\sf\large\leftmark}
 \rule[-.6ex]{0mm}{2mm} \hfill\sf\mdseries{\thepage}}
 \vskip 1pt \hrule}}}

\def\ps@myfancyplain{
 \def\@evenhead{\vbox{\hsize\textwidth%
 \rule[-.6ex]{0mm}{2mm} \hfill }
 \vskip 1pt \hrule
 \vskip\headsep
 \vskip\textheight
 \vskip1pc
 \hbox to \textwidth{\sf\mdseries\thepage 
 \rule[-6ex]{0mm}{2mm} \hfill }}
 \def\@oddhead{\vbox{\hsize\textwidth 
 \vskip 1pt\hrule
 \vskip\headsep
 \vskip\textheight
 \vskip2pc
 \hbox to \textwidth{\hfill\rule[.4ex]{1pc}{2.5pt}
 \sf\mdseries\thepage}
}}}

\def\ps@myemptyfun{\ps@plain}

\makeatother
\providecommand{\proofname}{Demostraci\'on.}
 {\par\noindent{\it Demostraci\'on. }\nopagebreak\normalsize}%
 {}

 {\par\noindent{\it #1. }\nopagebreak\normalsize}%
 {\hfill\linebreak[2]\hspace*{\fill}$\square$\\[-1pt]}
\makeatother

\def\sqbullet{\raise.2ex\hbox{\vrule width 3.5pt height 3.5pt}}

{}

{\em}

\newcounter{substep}
\def\thesubstep{\arabic{substep}}

{\em}

\newcounter{subsubstep}
\def\thesubsubstep{\arabic{subsubstep}}

{\em}

\numberwithin{figure}{section}

{}

\newtheoremstyle{mystyle}
  {}
  {}
  {\itshape}
  {}
  {\sf \bfseries}
  {}
{ }
  {\thmname{#1}\thmnumber{{\textcolor{blue}{\, \hspace{-1mm}#2.}}}\thmnote{ (#3)}}

\theoremstyle{mystyle}
\definecolor{royalblue(web)}{rgb}{0.25, 0.41, 0.88}
\hypersetup{
    colorlinks=true,
    linkcolor=royalblue(web),
    filecolor=magenta,      
    urlcolor=cyan,
    pdftitle={Overleaf Example},
    pdfpagemode=FullScreen,
    }

\titleformat{\section}
  {\normalfont\LARGE\bfseries}{\thesection.}{1em}{}
\titleformat{\subsection}
  {\normalfont\Large\bfseries}{\thesubsection.}{1em}{}
\titleformat{\subsubsection}
  {\normalfont\normalsize\itshape}{\thesubsubsection.}{1em}{}

\newtheorem{Teor}{Theorem}[section]

\newtheorem{Prop}[Teor]{Proposition}

\newtheorem{Lema}[Teor]{Lemma}

\newtheorem{Obse}[Teor]{Remark}

\newcommand{\norm}[1]{\left\lVert#1\right\rVert}

\newcommand{\R}{{\mathbb R}}

\newcommand{\dist}{\operatorname{dist}}

\newcommand{\supp}{\operatorname{supp}}

\newcommand{\id}{\operatorname{id}}

\newcommand{\mail}[1]{\small\href{mailto:#1}{#1}}
\newcommand{\SO}{\operatorname{SO}(n)}

\newcommand{\Ds}{D_{\rho,\mathrm{sym}}}
\newcommand{\Hz}{H_0^{\rho}(\Omega;\R^n)}
\newcommand{\Hr}{H^{\rho}(\R^n;\R^n)}
\newcommand{\dd}{\,dx}

\newenvironment{Abstract}
{
\begin{center}
\textbf{Abstract}\\
\vspace{0.25cm}
\begin{minipage}{14.5cm}}
{\footnotesize
\end{minipage}
\end{center}}

\begin{document}


	\begin{center}
		{\huge {\bfseries Linearization of nonlocal elasticity based on nonlocal gradients via $\Gamma$-convergence}\par}
		\vspace{1cm}
		
\begin{tabular}{l@{\hskip 2cm}l} 
	{\Large José C. Bellido}{\small\textsuperscript{1}} & {\Large Guillermo García-Sáez}{\small\textsuperscript{1}} \\
	\mail{josecarlos.bellido@uclm.es} & \mail{guillermo.garciasaez@uclm.es}
\end{tabular}
\vspace{5mm}

\textsc{\textsuperscript{1}ETSII, Departamento de Matem\'aticas,\\ Universidad de Castilla-La Mancha.} \\
		Campus Universitario s/n, 13071 Ciudad Real, Spain. \\ \vspace{5mm}
\end{center}
\begin{Abstract}
We provide a rigorous variational derivation of linear nonlocal elasticity from a general nonlinear elastic model based on nonlocal gradients. Under standard assumptions of nonlinear elasticity, we prove that the rescaled nonlinear energies $\Gamma$-converge, in the small-strain limit, to the corresponding quadratic linearized energy. The proof relies on a nonlocal geometric rigidity estimate and does not require any global upper growth condition on the stored-energy density. 
\end{Abstract}



\tableofcontents


\section{Introduction} The passage from nonlinear elasticity to linearized elasticity is a fundamental problem in the variational analysis of elastic materials. At the formal level, the linearized model is obtained by considering small perturbations of a reference configuration and expanding the stored-energy density around the identity. A rigorous justification, however, requires more than a pointwise Taylor expansion. In particular, one must control the nonlinear energies uniformly in the small-strain parameter and identify their variational limit on the appropriate function space. In the classical local setting, the rigorous derivation of linearized elasticity from finite elasticity can be formulated in terms of $\Gamma$-convergence. A central ingredient in this analysis is the geometric rigidity estimate of Friesecke, James, and Müller \cite{FJM}, which allows the nonlinear coercivity condition 
\[ W(F)\geq c\,\operatorname{dist}^{2}(F,SO(n)) \] 
to be converted into an estimate for the linearized strain. This approach was used by Dal Maso, Negri, and Percivale \cite{DMNP} to establish linearized elasticity as the $\Gamma$-limit of finite elasticity under natural assumptions on the stored-energy density. The purpose of the present paper is to develop an analogous rigorous linearization result for a broad class of nonlocal elastic models defined through nonlocal gradients. More precisely, given a radial kernel $\rho$, we consider the nonlocal gradient 
\[ D_\rho v(x) = \int_{\mathbb{R}^{n}} \frac{v(x)-v(y)}{|x-y|} \otimes \frac{x-y}{|x-y|} \rho(x-y)\,dy. \] 
This operator provides a unified first-order framework that includes, among other examples, the Riesz fractional gradient and truncated fractional gradients motivated by finite-horizon interactions in peridynamics (see \cite{ Bellido_NH2020, BellidoCuetoMoraCorral2021, BellidoCuetoMoraCorral2023,Bellido_NH2024,ShiehSpector2015,ShiehSpector2018} and the references therein). Its associated function spaces retain several essential features of classical Sobolev spaces, including density results, Poincaré-type inequalities, compact embeddings, and a suitable nonlocal Korn inequality. The reference \cite{BMS} addresses these structural properties for the spaces associated to nonlocal gradients for general kernels $\rho$. 


A formal derivation of nonlocal linear elasticity from a general hyperelastic model based on these nonlocal gradients was recently developed by the authors in \cite{BGS}. In that work, the linearized equilibrium equations were obtained by a formal expansion of the nonlinear constitutive law around the identity. The resulting model takes the form 
\begin{equation}\label{eq:limitpb} \begin{cases} -\operatorname{div}_{\rho} \bigl({C}[D_{\rho,\mathrm{sym}}v]\bigr)=f &\text{in }\Omega,\\ v=0 &\text{on }\Omega^{c}, \end{cases} \end{equation}
where \[ D_{\rho,\mathrm{sym}}v := \frac{D_\rho v+D_\rho v^{T}}{2} \] and ${C}=D^{2}W(I)$ is the elasticity tensor associated with the stored-energy density. It is worth pointing out that nonlocal linear model \eqref{eq:limitpb} fits into Eringen's theory of nonlocal elasticity \cite{BellidoCuetoMora2023,BGS,Eringen2002, EvgrafovBellido2019}. 

The aim of this paper is to provide a rigorous variational justification of the previous formal linearization, thus, this paper can be seen as a follow-up of \cite{BGS}. We consider a bounded reference configuration $\Omega\subset\mathbb{R}^{n}$ and deformations of the form \[ u=\operatorname{id}+hv, \] where $h>0$ is the small-strain parameter and $v\in \Hz$. The corresponding rescaled energies are 
\[ E_h(v) = \frac{1}{h^{2}} \int_{\mathbb{R}^{n}} W\bigl(I+hD_\rho v\bigr)\,dx - \int_{\Omega}f\cdot v\,dx. \] 
Under the standard assumptions of nonlinear elasticity, namely frame indifference, a stress-free reference configuration, and quadratic coercivity with respect to $SO(n)$, we prove that $(E_h)$ $\Gamma$-converges, with respect to the strong topology of $L^{2}(\mathbb{R}^{n};\mathbb{R}^{n})$, to the quadratic functional 
\[ E(v) = \frac12 \int_{\mathbb{R}^{n}} {C}[D_{\rho,\mathrm{sym}}v] :D_{\rho,\mathrm{sym}}v\,dx - \int_{\Omega}f\cdot v\,dx. \] 
The limit functional is precisely the energy associated with the linear nonlocal elasticity system derived formally in \cite{BGS}. The result is obtained without imposing a global upper growth condition on $W$. In particular, the energy density may exhibit superquadratic growth or take the value $+\infty$ outside a neighborhood of the identity. This is relevant in nonlinear elasticity, where one may wish to exclude matrices with nonpositive determinant or to incorporate singular behavior associated with extreme compression. Our result follows the lines of \cite{DMNP}, being the main ingredient the nonlocal version of the geometric rigidity given in Proposition \ref{prop:rigidity}. 

While this work was being completed, we became aware of the related contribution \cite{KreisbeckMoraSeifert}, where the linearization of a class of nonlocal hyperelastic energies is also rigorously established by means of $\Gamma$-convergence. Although the two works share a similar objective, their approaches are substantially different. The argument in \cite{KreisbeckMoraSeifert} relies on a local-to-nonlocal translation method: the identity 
\[ D_\rho u=\nabla(Q_\rho*u) ,\]
where $Q_\rho$ is the {\it translation} potential (see Section 2.2), is used to reduce the nonlocal problem to an available linearization result in classical elasticity. In contrast, our analysis is carried out directly within the nonlocal framework and is based on a nonlocal geometric rigidity estimate inspired by the celebrated rigidity theorem in \cite{FJM}. This result is obtained from the local one in \cite{FJM} by making use of the translation method as well. The two approaches therefore provide complementary rigorous justifications of the small-strain limit in nonlocal elasticity. We also mention the recent work \cite{Scott2026}, which pursues a similar objective by establishing rigidity estimates in nonlocal Gagliardo-type spaces.

The paper is organized as follows. In Section~2, we introduce the nonlocal gradient, its relationship with classical gradient, and the nonlocal Sobolev spaces used throughout the paper. We also recall the density, Poincaré, compactness, translation, and Korn-type results needed in the analysis. In Section~3, we define the rescaled nonlinear energies and state the main $\Gamma$-convergence theorem. Its proof is divided into the derivation of the nonlocal geometric rigidity estimate, the equicoercivity argument, the lower-bound inequality, and the construction of recovery sequences. Finally, in Section~4, we prove existence and uniqueness for the limit problem \eqref{eq:limitpb} (this was already shown in \cite{BGS} for general nonlocal gradients and in \cite{BellidoCuetoMora2023} for Riesz fractional gradients) and establish the convergence of minimum values and almost minimizing sequences.

\section{Preliminaries}\label{sec:kernel}

\subsection{Notation}
We fix $n\in \mathbb{N}$ the dimension of our ambient space $\R^n$ and we will denote by $\Omega\subset \R^n$ an open, bounded subset representing the body configuration. The notation for Sobolev $W^{1,p}$ and Lebesgue $L^p$ spaces is the standard one, as it is that of smooth
functions of compact support $C_c^\infty$. We will
indicate the domain of the functions, as in $L^p
(\Omega$); the target is indicated only if it is not $\R$.

We write $|x|=\left(\sum_{j=1}^n x_j^2\right)^{1/2}$ for the Euclidean norm of a vector $x=(x_1,\ldots,x_n)\in \R^n$, and $|A|$ for the Frobenius norm of a matrix $A\in \R^{m\times n}$. The identity matrix will be denoted as $I$ and the identity map as $\operatorname{id}.$ The subset of the symmetric matrices is denoted by $\R^{n\times n}_{sym}$, and the symmetric part of a matrix $M\in \R^{n\times n}$, will be denoted as $$M_{sym}=\frac{M+M^T}{2},$$ where $M^T$ is the transpose matrix of $M$. The inner product in $\R^n$ will be denoted by $\cdot$, while the respective one in $\R^{n\times n}$ will be denoted by $:$, i.e., given $A=(a_{ij})_{1\leq i,j\leq n}$ and $B=(b_{ij})_{1\leq i,j\leq n}$, $$A:B=\sum_{i,j=1}^n a_{ij}b_{ij}.$$
The tensor product of two vectors in $\R^n$ will be denoted by $\otimes$, and $*$ represents the usual operation of convolution.

The ball centered at $x\in \R^n$ and with radius $r>0$ is denoted by $B(x,r)=\{y\in \R^n: |x-y|<r\}$. The complementary of a set $E\subset \R^n$ is denoted by $E^c=\R^n\setminus{E}$, its closure by $\overline{E}$ and its boundary by $\partial E$. The distance between a point $x\in \R^n$ and an open set $E\subset\R^n$ is denoted by $d(x,E)$.

Our convention for the Fourier transform of functions $f\in L^1(\R^n)$ is $$\widehat{f}(\xi)=\int_{\R^n}f(x)e^{-2i\pi x\cdot \xi}\,dx,\,\xi\in \R^n.$$ This definition is extended by continuity and duality to other function and distribution spaces as usually in function spaces theory. We will sometimes use the alternative notation $\mathcal{F}(f)$ for $\widehat{f}$. For more details on Fourier analysis we refer to the classical text \cite{Grafakos2008}.

Regarding radial functions we have the following definitions:
\begin{itemize}
    \item A function $f:\R^n\to \R$ is \textit{radial} if there exists $\overline{f}:[0,\infty)\to \R$ such that $\overline{f}(|x|)=f(x)$ for every $x\in \R^n$. The function $\overline{f}$ is called the \textit{radial representation} of $f$.
    \item A radial function $f:\R^n\to \R$ is \textit{radially decreasing} if its radial representation is a
decreasing function.
\item A function $g:\R^n\to \R^n$ is \textit{vector radial} if there exists a radial function $\overline{g}:[0,\infty)\to \R$  such that
$g(x)=\overline{g}(|x|)x$ for every $x\in \R^n$.
\end{itemize}

For real valued functions, we use the monotonicity properties of being increasing and decreasing in
the non-strict sense. A function $f :\R \to\R$ is called \textit{almost decreasing} if there exists a positive constant $C$ such
that $f(t)\geq Cf(s)$ for every $t\leq s$, and an analogous definition for \textit{almost increasing}.

A function $f:\Omega\times \R^m\times \R^{m\times n}$ is called a \textit{Carath\'eodory integrand} if $f(\cdot, z,A)$ is measurable for all $(z,A)\in \R^m\times \R^{m\times n}$ and $f(x,\cdot,\cdot)$ is continuous a.e. $x\in \Omega$.

\subsection{Nonlocal gradients}

Throughout the entire manuscript, we will assume that the kernel, $\rho$, in the definition of nonlocal gradients is a function that verifies:
\begin{equation}\label{eqn:H0} \tag{H0}\begin{cases} \text{$\rho:\R^n\setminus\{0\} \to [0,\infty)$ is radial}, \\
\text{$\rho \in L^1_{\text{loc}}(\R^n)$ with $\displaystyle \int_{\R^n} \min\{1,|x|^{-1}\}\rho(x)\,dx < \infty$,} \\
\text{$\inf_{\overline{B(0,\varepsilon)}}\rho >0$ for some $\varepsilon>0$.} \end{cases}\end{equation}
Hypotheses, as (\ref{eqn:H0}), are somehow the minimal requirement for the kernel $\rho$ in the radial case and have frequently appeared in the related literature. Observe that, since by \eqref{eqn:H0} the kernel $\rho$ is radial, there exists a function $\overline{\rho}:(0,+\infty)\to [0,+\infty)$ such that $\rho(x)=\overline{\rho}(|x|)$. With a slight abuse of notation, and without risk of ambiguity, we shall henceforth denote $\overline{\rho}$ simply by $\rho$. We also assume, as in the
linearization setting of \cite{BGS}, that
\begin{equation}\label{eq:normalization}
 \supp\rho\text{ is compact, and}\qquad
 \int_{\R^n}\rho(x)\dd=n.
\end{equation}

For $v\in C_c^\infty(\R^n;\R^n)$, the nonlocal gradient based on the kernel $\rho$ and its symmetrization are defined as,
\begin{equation}
 D_\rho v(x)=\int_{\R^n}
 \frac{v(x)-v(y)}{|x-y|}\otimes
 \frac{x-y}{|x-y|}\rho(x-y)\,dy,
 \qquad D_{\rho,\mathrm{sym}}v=(D_\rho v)_{\mathrm{sym}}.
\end{equation}

The associated potential is defined by 
\begin{equation}\label{eq:Q} Q_\rho(x) = \int_{\lvert x\rvert}^{\infty}\frac{\rho(r)}{r}\,dr, \qquad x\in\mathbb{R}^n\setminus\{0\}. \end{equation} 
It is radial, nonnegative and compactly supported. Moreover, assuming $\rho\in L^1(\R^n)$, by Tonelli's theorem and polar coordinates, 
 \begin{equation}\label{eq:mass} \begin{aligned}\|Q_\rho\|_{L^1(\mathbb{R}^n)} &= \lvert\mathbb{S}^{n-1}\rvert \int_0^\infty r^{n-1} \left( \int_r^\infty \frac{\rho(s)}{s}\,ds \right)dr \\ &= \frac{\lvert\mathbb{S}^{n-1}\rvert}{n} \int_0^\infty \rho(s)s^{n-1}\,ds \\ &= \frac1n\int_{\mathbb{R}^n}\rho(x)\,dx = 1. \end{aligned}\end{equation}  
 In particular, \(Q_\rho\in L^1(\mathbb{R}^n)\).

A major result on the theory of nonlocal gradients are the translation property by the associated potential, which allows us to see nonlocal gradients as translation of classical gradients. In particular, by \cite[Proposition 2.6]{BMS} for a function $v\in C^\infty_c(\R^n;\R^n)$, it holds,
\begin{equation}\label{eq:forward}
 D_\rho v=Q_\rho*Dv=D(Q_\rho*v).
\end{equation}
The operator $\operatorname{div}_\rho$ is the negative
distributional adjoint of $D_\rho$ and can be defined for $v\in C_c^\infty(\R^n,\R^{n\times n})$ as  $$\operatorname{div}_\rho v(x):=\int_{\R^n}\frac{v(x)-v(y)}{|x-y|} \frac{x-y}{|x-y|}\rho(x-y)\,dy,\,x\in \R^n.$$ As in the local case, we have an integration by parts
identity \cite[Proposition 3.2]{BMS}. 
\begin{Lema}[Integration by parts]\label{IBP}
    Let $u\in C_c^\infty(\R^n)$ and $v\in C_c^\infty(\R^n;\R^n)$. Then, $$\int_{\R^n}D_\rho u\cdot v\,dx=-\int_{\R^n}u\operatorname{div}_\rho v\,dx.$$
\end{Lema}

\subsection{Nonlocal Sobolev spaces}
Analogously as the classical and fractional cases, the \textit{nonlocal Sobolev spaces} based on $D_\rho$ were introduced in \cite{BMS},
\begin{align*}
 H^{\rho,p}(\R^n;\R^n)=\{v\in L^p(\R^n;\R^n):
       D_\rho v\in L^p(\R^n;\R^{n\times n})\},
\end{align*} and the closed proper subspace $$ H^{\rho,p}_0(\Omega;\R^n)=\{v\in H^{\rho,p}(\R^n;\R^n):v=0\text{ a.e. on }\Omega^c\},$$ for $1\leq p\leq \infty$ and $\Omega\subset\R^n$ and open subset. Under the norm $$\|v\|_{\rho,p}:=\begin{cases}\left(\|v\|_p^p+\|D_\rho v\|_p^p\right)^\frac{1}{p},\,1\leq p<\infty,\\
\operatorname{max}\{\norm{u}_\infty,\norm{D_\rho u}_\infty\},\,p=\infty,\end{cases},$$
the spaces are complete. In the hilbertian case $p=2$, we will denote the spaces as $\Hr$ and $\Hz$, respectively. 

For $1\leq p<\infty$, it was proven in \cite[Theorem 3.9]{BMS} the density results:
\begin{equation}\label{eq:density}
 \overline{C_c^\infty(\R^n;\R^n)}^{\,\|\cdot\|_{\rho,p}}=H^{\rho,p}(\R^n;\R^n),\qquad\overline{C_c^\infty(\Omega;\R^n)}^{\,\|\cdot\|_{\rho,p}}=H_0^{\rho,p}(\Omega;\R^n),
\end{equation} under the assumption of $\Omega$ being Lipschitz.

We recall that we can restrict ourselves to kernels with compact support without loss of generality since the main properties of the function spaces $H^{\rho,p}(\R^n)$ and $H_0^{\rho,p}(\Omega)$ are encoded on the behavior of $\rho$ near
zero by \cite[Proposition 3.10]{BMS}.

In order to ensure the fundamental structural properties of the spaces $H^{\rho,p}$, such as Poincar\'e-type inequalities and compact embeddings, it is necessary to impose additional assumptions on the kernels $\rho$.

Let $\varepsilon$ be as in \eqref{eqn:H0}. We will employ the following conditions:
\begin{enumerate}
    \item[(H1)]\label{eqn:H1} There exists $\nu>0$ such that the functions $f_\rho:(0,\infty)\to \R$, defined by $f_\rho(r)= r^{n-2}{\rho}(r)$, and $g(r)= r^\nu f_\rho(r)$ are decreasing on $(0,\varepsilon)$;
    \item[(H2)]\label{eqn:H2} The function $f_\rho$ is smooth away from the origin, and for every positive integer $k$ there exists a constant $C(k)>0$ such that

\[
    \left|\frac{d^kf_\rho(r)}{dr^{k}}\right|\leq C(k)\,r^{-k}f_\rho(r), \quad r\in (0,\varepsilon);
    \]

    \item[(H3)]\label{eqn:H3} There exists $s\in(0,1)$ such that the function $r\mapsto r^{n+s-1}{\rho}(r)$ is almost decreasing on $(0,\varepsilon)$;
    \item[(H4)]\label{eqn:H4} There exists $t\in(0,1)$ such that the function $r\mapsto r^{n+t-1}{\rho}(r)$ is almost increasing on $(0,\varepsilon)$.
\end{enumerate}
These hypotheses are natural, and we refer to Sections~4 and~5 of \cite{BMS} for a detailed discussion.

In \cite{BMS, CuetoHidde2025}, examples of admissible kernels are given, the most paradigmatic one 
	being $$\rho^s(x):=c_{n,s}\frac{1}{|x|^{n+s-1}},\,s\in (0,1),$$ where $c_{n,s}$ is a suitable normalization constant. The nonlocal gradient associated with this kernel is the well-known Riesz fractional gradient $D^s=:D_{\rho^s}$, which is defined for functions $u\in C_c^\infty(\R^n)$ as $$D^su(x)=c_{n,s}\int_{\R^n}\frac{u(x)-u(y)}{|x-y|^{n+s}}\frac{x-y}{|x-y|}\,dy=D(I_{1-s}u),$$ where $I_s$ is the Riesz potential.
	The Riesz fractional gradient was introduced by Shieh and Spector in \cite{ShiehSpector2015,ShiehSpector2018}, and since then has been widely studied.
	
	Another important example concerns a truncated version of the Riesz fractional gradient introduced in \cite{BellidoCuetoMora2023}, which is suitable for bounded domains. This formulation is motivated by concepts from state-based peridynamics, as its integration domain is restricted to a ball of radius $\delta>0$, corresponding to the interaction horizon among particles in the peridynamic framework. Let $\delta>0$, and a non-negative radial cut-off
	function $w_\delta$ such that $w_\delta(0)>0$, $w_\delta\in C_c^\infty\left(B(0,\delta)\right)$ and $w_\delta/|\cdot|^{1+s}$ is radially decreasing.
	Now, taking the kernel 
	$$\rho_\delta^s(x):=\frac{w_\delta(x)}{|x|^{n+s-1}},$$ 
	yields that $D_{\rho_\delta^s}$ is, up to a normalization constant, equal to $D_\delta^s$, the nonlocal gradient introduced in \cite{BellidoCuetoMora2023}.

Under these hypothesis, Poincar\'e and compactness results were proven for general kernels with suitable behaviour around the origin.
\begin{Teor}\label{eq:poincare}
    Let $p\in (1,\infty)$, $\Omega\subset \R^n$ an open bounded set and $\rho$ satisfying (H0)-(H2) and having compact support. Then, the following two statements hold:
    \begin{enumerate}
        \item (Poincar\'e inequality) If $$\liminf_{t\to 0^+}t^{n-1}\overline{\rho}(t)>0,$$ then there is a positive constant $C$ such that $$\norm{u}_{L^p(\Omega)}\leq C\norm{D_\rho u}_{L^p(\R^n;\R^n)}\,$$
        for all $u\in H_0^{\rho,p}(\Omega)$
        \item (Compact embedding) If $$\lim_{t\to 0^+}t^{n-1}\overline{\rho}(t)=\infty,$$ then $H^{\rho,p}_0(\Omega)$ is compactly embedded into $L^p(\R^n)$.
    \end{enumerate}
\end{Teor}
 Assuming (\ref{eqn:H0})--(H4) for $\rho$, we obtain the following lemma, which establishes the existence of some translation operators from the nonlocal framework to the local one and vice versa \cite[Lemma~2.12]{CuetoHidde2025}, extending the representation formula (\ref{eq:forward}) to functions on $\Hr$.
	\begin{Lema}[Translation operators]\label{eq:forward-space}
		Let $\rho$ satisfy (\ref{eqn:H0})--$(H4)$ with compact support. The linear map $\mathcal{Q}_\rho$ defined as $\mathcal{Q}_\rho u:=Q_\rho *u$, $u\in C_c^\infty(\R^n)$, extends to a linear bounded map from $H^{\rho,p}(\R^n)\to W^{1,p}(\R^n)$. Furthermore, it holds for all $u\in H^{\rho,p}(\R^n)$ that $$D_\rho u=D\left(\mathcal{Q}_\rho u\right).$$ The linear map $\mathcal{P}_\rho$ defined as $$v\mapsto\mathcal{P}_\rho  v:=\mathcal{F}^{-1}\left(\widehat{v}/\widehat{Q}_\rho\right),\,v\in \mathcal{S}(\R^n),$$ extends to a bounded operator from $W^{1,p}(\R^n)$ to $H^{\rho,p}(\R^n)$. Moreover, $(\mathcal{P}_\rho)^{-1}=\mathcal{Q}_\rho$, and for every $v\in W^{1,p}(\R^n)$, $$D v=D_\rho(\mathcal{P}_\rho v).$$
	\end{Lema}
Using this operators, in \cite[Proposition 4.2]{BGS} was proven the following Korn's inequality on the nonlocal setting.
\begin{Prop}[Nonlocal Korn's inequality]\label{eq:korn}
    Let $1<p<\infty$. Then, there exists $c_p>0$ such that for all $u\in H^{\rho,p}(\R^n;\R^n)$, $$\norm{D_{\rho,sym}u}_{L^p(\R^n)}\geq c_p\norm{D_\rho u}_{L^p(\R^n)}.$$
\end{Prop}
\section{A general model via linearization}\label{sec:model}
The key assumptions of this paper are the usual assumptions of nonlinear elasticity theory for the stored energy following \cite[pp. 1466]{FJM} . In particular, 
\label{sec:W}
we let $W:\R^{n\times n}\to[0,+\infty]$ be Borel representing the stored-energy density
and we assume that it satisfies:
\begin{description}
\item[(W1)] $W$ is finite and $C^2$ on an open neighborhood of $I$.
\item[(W2)] $W$ is stress free and frame indifferent, i.e., $W(I)=0$ and $W(RF)=W(F)$ for every rotation $R\in\SO$
and every $F\in\R^{n\times n}$.
\item[(W3)] There exists $c>0$ such that
\begin{equation}\label{eq:Wcoercivity}
 W(F)\ge c\,\dist^2(F,\SO),\qquad F\in\R^{n\times n}.
\end{equation}
\end{description}
We do not impose a global upper quadratic growth condition.
In particular, $W$ may take the value $+\infty$ outside its
neighborhood of regularity, including at matrices with nonpositive
determinant. Conditions ensuring existence of minimizers of each
nonlinear energy are not needed for the convergence theorem.

Following \cite[Section 3]{BGS}, we now introduce the following hyperelastic nonlinear energy associated with a given deformation mapping 
$u \in H^{\rho,p}_{\operatorname{id}}(\Omega;\R^n)$:
\[
E[u]=\int_\Omega W\bigl(x,D_\rho u\bigr)\,dx - \int_\Omega f \cdot u\,dx,
\]
where the integrand $W:\Omega \times \R^{n\times n} \to \R \cup \{+\infty\}$ is a Carath\'eodory function representing the stored elastic energy density on the previous hypothesis, and $f:\Omega \to \R^n$ denotes a distributed body force acting on the domain. For simplicity, we impose Dirichlet boundary conditions on the complementary domain, namely $u|_{\Omega^c}=\operatorname{id}$.

The Dirichlet problem associated to the minimization problem for $E$ is $$\begin{cases}
    -\operatorname{div}_\rho W_2(x,D_\rho u)&=f,\,\text{in $\Omega$},\\
    \hfill u&=\operatorname{id},\,\text{on $\Omega^c$}
\end{cases}$$
where $W_2(x,F)=:T_R(x,F)$ is the Piola-Kirchhoff stress tensor, which depends on nonlinear operators as the deformation jacobian. In order to obtain simpler linear models to work with, in \cite{BGS} it was performed a formal linearization process which yields the strong form of the equations of general nonlocal linear elasticity: \begin{equation}\label{eqn:P}\begin{cases}
    -\operatorname{div}_\rho \left(C[D_{\rho,sym}v]\right)&=f,\,\text{in $\Omega$},\\
    \hfill v&=0,\,\text{on $\Omega^c$}.
\end{cases}\end{equation}
Multiplying both sides by a function $w\in H^{\rho}_0(\Omega;\R^n)$, and using again the duality relationships of Lemma \ref{IBP}, we get that the weak form of the equations is the following: find $v\in H^\rho_0(\Omega;\R^n)$ such that 
\begin{equation*}\int_{\R^n} C[D_{\rho,sym}v]:D_{\rho,sym}w\,dx=\int_{\R^n} f\cdot w\,dx,\,\forall w\in H^{\rho}_0(\Omega;\R^n).\end{equation*}
Here, $C$ is the elasticity tensor, obtained as the gradient of $T_R$ evaluated on $I$. Since $T_R(F)=W_2(F)=DW(F)$, we have that the components of the tensor are 
$$c_{ijkl}=\frac{\partial^2 W}{\partial F_{ij}\partial F_{kl}}(I),\quad 1\leq i,j,k,l\leq n.$$
$C$ is a fourth order tensor acting over matrices $M$, we will denote the action of $C$ over $M$ as $C[M]$ instead of $CM$ to differentiate it from the usual product of matrices. One of the main properties of the elasticity tensor is that $C[M]=C[M_{sym}]\in \R^{n\times n}_{sym}$, hence we can consider $C$ as a map from $\R_{sym}^{n\times n}\to \R^{n\times n}_{sym}$. As in the classical case, $C[D_\rho] v=C[D_{\rho,sym} v]$, where $D_{\rho,sym}v=(D_\rho v)_{sym}$. 

We now state a technical lemma that will be useful later. 
\begin{Lema}\label{lem:quadratic}
The elasticity tensor is self-adjoint
for the Frobenius product, and there exist $c_1,c_2>0$ such that
\begin{equation}\label{eq:quadratic-bounds}
 c_1|M_{\mathrm{sym}}|^2
 \le\frac12C[M]:M
 \le c_2|M_{\mathrm{sym}}|^2.
\end{equation}
In particular, the quadratic form $M\mapsto\tfrac12C[M]:M$ is convex. Moreover, there exists a nondecreasing function $\omega$, with
$\omega(r)\to0$ as $r\downarrow0$, such that for $|M|$ small,
\begin{equation}\label{eq:Taylor}
 \left|W(I+M)-\frac12C[M]:M\right|
 \le\omega(|M|)|M|^2.
\end{equation}
\end{Lema}
\begin{proof}
Since $W\ge0$ and $W(I)=0$, we have $DW(I)=0$.
Symmetry of the Hessian gives
\[
C[M]:N=C[N]:M,
\]
so $C$ is self-adjoint.

Let $A^T=-A$ and $N\in\mathbb R^{n\times n}$.
By frame indifference,
\[
W\bigl(e^{tA}(I+\varepsilon N)\bigr)=W(I+\varepsilon N).
\]
Differentiating first in $t$ and then in $\varepsilon$ at
$(t,\varepsilon)=(0,0)$ yields
\[
C[N]:A+DW(I):AN=0.
\]
Thus $C[N]:A=0$ for every skew-symmetric $A$, so $C[N]$ is symmetric.
Self-adjointness also gives $C[A]=0$. Consequently,
\[
C[M]=C[M_{\mathrm{sym}}],
\qquad
C[M]:M=C[M_{\mathrm{sym}}]:M_{\mathrm{sym}}.
\]

For a symmetric matrix $S$, the matrix $I+hS$ is positive definite
when $|h|$ is sufficiently small. Its polar decomposition gives
\[
\operatorname{dist}(I+hS,SO(n))=|hS|.
\]
Writing $c_W>0$ for the constant in \textnormal{(W3)}, Taylor expansion
and coercivity imply
\[
\frac12 C[S]:S
=\lim_{h\to0}\frac{W(I+hS)}{h^2}
\ge c_W|S|^2.
\]
The upper bound follows from
$|C[S]:S|\le\|C\||S|^2$.
Hence \eqref{eq:quadratic-bounds} holds with
$c_1=c_W$ and $c_2=\frac12\|C\|$.

Choose $r_0>0$ such that $W$ is $C^2$ on a neighborhood of
$\{I+N:|N|\le r_0\}$, and set
\[
\omega(r):=\frac12
\sup_{|N|\le r}\|D^2W(I+N)-D^2W(I)\|.
\]
This function is nondecreasing and tends to zero at the origin.
Taylor's integral formula gives
\[
W(I+M)-\frac12 C[M]:M
=
\int_0^1(1-t)
\bigl(D^2W(I+tM)-D^2W(I)\bigr)[M,M]\,dt,
\]
which proves \eqref{eq:Taylor}.
\end{proof}

\subsection{Rescaled energies and the main result}\label{sec:theorem}
The globally defined reference deformation satisfies
\begin{equation}\label{eq:affine}
 D_\rho\id=\int_{\R^n}\frac{z\otimes z}{|z|^2}\rho(z)\,dz=I
\end{equation}
by radiality and \eqref{eq:normalization}. We therefore consider
\begin{equation}\label{eq:deformation}
 u=\id+hv,\qquad v\in\Hz,\qquad D_\rho u=I+hD_\rho v.
\end{equation}
This is an affine admissible class: $u=\id$ on $\Omega^c$.

Fix $f\in L^2(\Omega;\R^n)$. The rescaled problem on
$L^2(\R^n;\R^n)$ is
\begin{equation}\label{eq:Eh}
 E_h(v)=
 \begin{cases}
 \displaystyle\frac1{h^2}\int_{\R^n}W(I+hD_\rho v)\dd
              -\int_\Omega f\cdot v\dd,&v\in\Hz,\\
 +\infty,&v\notin\Hz.
 \end{cases}
\end{equation}
As we have seen, the linearized energy is given by
\begin{equation}\label{eq:E}
 E(v)=
 \begin{cases}
 \displaystyle\frac12\int_{\R^n}C[D_{\rho,\mathrm{sym}}v]:
                         D_{\rho,\mathrm{sym}}v\dd
              -\int_\Omega f\cdot v\dd,&v\in\Hz,\\
 +\infty,&v\notin\Hz.
 \end{cases}
\end{equation}
We now state our main result:
\begin{Teor}\label{thm:main}
Under (H0)--(H4), \eqref{eq:normalization}, the domain assumptions in
Section \ref{sec:kernel}, and (W1)--(W3),
\begin{equation}\label{eq:Gamma}
 \Gamma\hbox{-}\lim_{h\downarrow0}E_h=E
 \quad\text{in the strong topology of }L^2(\R^n;\R^n).
\end{equation}
More precisely:
\begin{enumerate}
\item There exist $c_1,c_2>0$, independent of $h>0$, such that
\begin{equation}\label{eq:coercivity}
 E_h(v)\ge c_1\|v\|_{\rho,2}^2-c_2\|f\|_{L^2(\Omega)}^2,
 \qquad v\in\Hz.
\end{equation}
Consequently, the family of functionals $\{ E_h\}$ is equicoercive in $H^\rho_0(\Omega;\R^n)$. 
\item If $h_j\downarrow0$ and $v_j\to v$ in $L^2$, then
$E(v)\le\liminf_jE_{h_j}(v_j)$.
\item For every $v\in H^\rho_0(\Omega;\mathbb{R}^n)$ and every sequence $h_j\downarrow 0$, there exists a sequence $v_j\in C_c^\infty(\Omega;\mathbb{R}^n)$ such that 
\[ v_j\to v \quad\text{in }H^\rho(\mathbb{R}^n;\mathbb{R}^n) \quad \text{and}\quad E_{h_j}(v_j)\to E(v). \]
\end{enumerate}
\end{Teor}
We will prove each item separately in the following sections. 

\subsection{Geometric rigidity}\label{sec:rigidity}
We use the following classical result
\cite[Theorem 3.1]{FJM}, which states that for every ball $B\subset\R^n$, there exists
$c_B>0$ such that for every $w\in H^1(B;\R^n)$ there is a rotation $R\in SO(n)$ with
\begin{equation}\label{eq:FJM-classical}
 \|Dw-R\|_{L^2(B;\R^{n\times n})}
 \le c_B\|\operatorname{dist}(Dw,SO(n))\|_{L^2(B)}.
\end{equation}
The following proposition translates this local result to our nonlocal framework.
\begin{Prop}\label{prop:rigidity}
There exists $c>0$, depending on $n,\Omega,\rho$ and independent of
$h>0$, such that
\begin{equation}\label{eq:nonlinear-rigidity}
 \|D_\rho v\|_{L^2(\R^n;\R^{n\times n})}
 \le\frac{c}{h}
 \|\operatorname{dist}(I+hD_\rho v,SO(n))\|_{L^2(\R^n)}
 \qquad \text{for all }v\in\Hz.
\end{equation}
\end{Prop}
\begin{proof}
Let $r$ such that $$\overline{\Omega}+\operatorname{supp}Q_\rho\subset B(0,r),$$ and fix
\[
 B=B(0,r+2),\qquad
 U=B(0,r+2)\setminus\overline{B(0,r+1)}.
\]
Then $|U|>0$. By \eqref{eq:forward-space}, the map
\[
 x\longmapsto x+h(Q_\rho*v)(x)
\]
belongs to $H^1(B;\R^n)$, and its classical gradient is
$I+hD_\rho v$. Applying \eqref{eq:FJM-classical} to this map
produces a single matrix $R\in SO(n)$ such that
\begin{equation}\label{eq:FJM}
 \|I+hD_\rho v-R\|_{L^2(B;\R^{n\times n})}
 \le c_B\|\operatorname{dist}(I+hD_\rho v,SO(n))\|_{L^2(B)}.
\end{equation}
The rotation may depend on $h$ and $v$, but the constant $c_B$
does not.

We first control this rotation using the exterior annulus.
Note $D_\rho v=0$ a.e. on $U$.
Since $U\subset B$ and the squared integrand is nonnegative,
\begin{align}
 \|I+hD_\rho v-R\|_{L^2(B)}^2
 &=\int_B|I+hD_\rho v(x)-R|^2\dd\notag\\
 &\ge\int_U|I+hD_\rho v(x)-R|^2\dd\notag\\
 &=\int_U|I-R|^2\dd
 =|U|\,|I-R|^2.\label{eq:restriction}
\end{align}
Equivalently, the $L^2$-norm of the restriction of the field to
$U$ is no larger than its $L^2$-norm on $B$.
Taking square roots and using \eqref{eq:FJM} gives
\begin{equation}\label{eq:rotation}
 |I-R|\le\frac{c_B}{|U|^{1/2}}
 \|\operatorname{dist}(I+hD_\rho v,SO(n))\|_{L^2(B)}.
\end{equation}
The triangle inequality in $L^2(B;\R^{n\times n})$ yields
\begin{align}
 h\|D_\rho v\|_{L^2(B;\R^{n\times n})}
 &=\|(I+hD_\rho v-R)+(R-I)\|_{L^2(B;\R^{n\times n})}\notag\\
 &\le\|I+hD_\rho v-R\|_{L^2(B;\R^{n\times n})}
       +\|R-I\|_{L^2(B;\R^{n\times n})}\notag\\
 &=\|I+hD_\rho v-R\|_{L^2(B;\R^{n\times n})}
       +|B|^{1/2}|R-I|.\label{eq:triangle}
\end{align}
Substituting \eqref{eq:FJM} and \eqref{eq:rotation} separately into
the two terms on the right-hand side gives
\begin{equation}\label{eq:anchored}
 h\|D_\rho v\|_{L^2(B;\R^{n\times n})}
 \le c_B\left(1+\frac{|B|^{1/2}}{|U|^{1/2}}\right)
 \|\operatorname{dist}(I+hD_\rho v,SO(n))\|_{L^2(B;\R^{n\times n})}.
\end{equation}
Both fields vanish outside $B$, so their norms on $B$ equal their
norms on $\R^n$. Dividing by $h>0$ proves
\eqref{eq:nonlinear-rigidity}. The sets $B$ and $U$ are fixed,
hence the constant is independent of $h$ and $v$.
\end{proof}

\subsection{Equicoercivity}\label{sec:equicoercivity}
\begin{proof}[Proof of Theorem \ref{thm:main}(1)]
Let $h>0$ and $v\in H_0^\rho(\Omega;\mathbb R^n)$.
If $E_h(v)=+\infty$, the coercivity inequality is immediate.
We therefore assume that $E_h(v)<+\infty$.

Let $c_{\mathrm{rig}}>0$ be the constant in
Proposition \ref{prop:rigidity}, so that
\[
 \|D_\rho v\|_{L^2(\mathbb R^n)}
 \le \frac{c_{\mathrm{rig}}}{h}
 \|\operatorname{dist}(I+hD_\rho v,SO(n))\|_{L^2(\mathbb R^n)}.
\]
Squaring and rearranging gives
\[
 \frac{1}{h^2}
 \int_{\mathbb R^n}
 \operatorname{dist}^2(I+hD_\rho v,SO(n))\,dx
 \ge
 \frac{1}{c_{\mathrm{rig}}^2}
 \|D_\rho v\|_{L^2(\mathbb R^n)}^2.
\]
Let $c_W>0$ be the constant in assumption (W3), namely
\[
 W(F)\ge c_W\operatorname{dist}^2(F,SO(n)).
\]
Combining the preceding two inequalities yields
\begin{align}
 \frac{1}{h^2}
 \int_{\mathbb R^n}W(I+hD_\rho v)\,dx
 &\ge
 \frac{c_W}{h^2}
 \int_{\mathbb R^n}
 \operatorname{dist}^2(I+hD_\rho v,SO(n))\,dx
 \notag\\
 &\ge
 \frac{c_W}{c_{\mathrm{rig}}^2}
 \|D_\rho v\|_{L^2(\mathbb R^n)}^2.
 \label{eq:elastic-control}
\end{align}
Set
\[
 c_0:=\frac{c_W}{c_{\mathrm{rig}}^2}>0.
\]

Let $c_P>0$ be the constant
in the Poincar\'e inequality \eqref{eq:poincare}:
\[
 \|v\|_{L^2(\mathbb R^n)}
 \le c_P\|D_\rho v\|_{L^2(\mathbb R^n)}.
\]
Since $v=0$ almost everywhere on $\Omega^c$, the
Cauchy--Schwarz inequality gives
\begin{align*}
 \left|\int_\Omega f\cdot v\,dx\right|
 &\le
 \|f\|_{L^2(\Omega)}
 \|v\|_{L^2(\Omega)}\\
 &=
 \|f\|_{L^2(\Omega)}
 \|v\|_{L^2(\mathbb R^n)}\\
 &\le
 c_P\|f\|_{L^2(\Omega)}
 \|D_\rho v\|_{L^2(\mathbb R^n)}.
\end{align*}
Consequently, by the definition of $E_h$ and
\eqref{eq:elastic-control},
\[
 E_h(v)
 \ge
 c_0\|D_\rho v\|_{L^2(\mathbb R^n)}^2
 -
 c_P\|f\|_{L^2(\Omega)}
 \|D_\rho v\|_{L^2(\mathbb R^n)}.
\]

By Young's inequality in the form
\[
 ab\le \frac{\eta}{2}a^2+\frac{1}{2\eta}b^2,
 \qquad a,b\ge0,\quad \eta>0,
\]
with
\[
 a=\|D_\rho v\|_{L^2(\mathbb R^n)},\qquad
 b=c_P\|f\|_{L^2(\Omega)},\qquad
 \eta=c_0,
\]
we have that
\[
 c_P\|f\|_{L^2(\Omega)}
 \|D_\rho v\|_{L^2(\mathbb R^n)}
 \le
 \frac{c_0}{2}\|D_\rho v\|_{L^2(\mathbb R^n)}^2
 +
 \frac{c_P^2}{2c_0}\|f\|_{L^2(\Omega)}^2.
\]
Substituting this estimate, we obtain
\[
 E_h(v)
 \ge
 \frac{c_0}{2}\|D_\rho v\|_{L^2(\mathbb R^n)}^2
 -
 \frac{c_P^2}{2c_0}\|f\|_{L^2(\Omega)}^2
 \ge
 \frac{c_0}{2(1+c_P)^2}\|v\|_{\rho,2}^2
 -
 \frac{c_P^2}{2c_0}\|f\|_{L^2(\Omega)}^2.
\]
This proves \eqref{eq:coercivity}, with
\[
 c_1=\frac{c_0}{2(1+c_P)^2},
 \qquad
 c_2=\frac{c_P^2}{2c_0}.
\]
Both constants depend only on the fixed data
$n,\Omega,\rho,W$, and are independent of $h$, $v$ and $f$.

We now prove equicoercivity. In particular, if $h_j>0$ and
$\sup_j E_{h_j}(v_j)<+\infty$, then $(v_j)$ is bounded in
$H_0^\rho(\Omega;\mathbb R^n)$. Weak compactness in this
Hilbert space and the compact embedding into $L^2$ provide,
passing by a nonrelabeled subsequence, a function
$v\in H_0^\rho(\Omega;\mathbb R^n)$ such that
\begin{equation}\label{eq:compactness}
 \begin{gathered}
 v_j\rightharpoonup v
 \quad\text{in }H_0^\rho(\Omega;\mathbb R^n),\\
 v_j\longrightarrow v
 \quad\text{in }L^2(\mathbb R^n;\mathbb R^n).
 \end{gathered}
\end{equation}
Since $D_\rho:H_0^\rho(\Omega;\mathbb R^n)
\to L^2(\mathbb R^n;\mathbb R^{n\times n})$ is bounded
and linear, we also have
\[
 D_\rho v_j\rightharpoonup D_\rho v
 \quad\text{in }L^2(\mathbb R^n;\mathbb R^{n\times n}).
\]
This completes the proof.
\end{proof}

\subsection{The lower bound}\label{sec:liminf}

The following lemma proves the liminf inequality for our energy functional and follows from the cutoff and weak lower semicontinuity argument used in the proof of \cite[Theorem 6.1(i)]{FJM}, we include the proof for the sake of completeness and to keep the manuscript self-contained. 

\begin{Lema}\label{lem:liminf}
Let $B$ be a bounded measurable set, $h_j\downarrow0$, and
$M_j\rightharpoonup M$ in $L^2(B;\R^{n\times n})$. Under
(W1)--(W3),
\begin{equation}\label{eq:matrix-liminf}
 \liminf_{j\to\infty}\frac1{h_j^2}\int_B W(I+h_jM_j)\dd
       \ge\frac12\int_B C[M]:M\dd.
\end{equation}
\end{Lema}
\begin{proof}
Set $b=\sup_j\|M_j\|_2^2<\infty$ and define
\[
 A_j=\{x\in B:|M_j(x)|\le h_j^{-1/2}\},\qquad
 N_j=\mathbf1_{A_j}M_j.
\]
Chebyshev's inequality gives $|B\setminus A_j|\le bh_j\to0$.
For each fixed $\Phi\in L^2(B;\R^{n\times n})$,
\[
 \left|\int_B(M_j-N_j):\Phi\dd\right|
 \le b^{1/2}\|\mathbf1_{B\setminus A_j}\Phi\|_2\longrightarrow0.
\]
The last step follows from absolute continuity of the integral of
$|\Phi|^2$. Consequently $N_j\rightharpoonup M$ in $L^2$.
On $A_j$ we have $|h_jM_j|\le\sqrt{h_j}$, so
\eqref{eq:Taylor} gives
\[
 \frac1{h_j^2}W(I+h_jM_j)
 \ge\frac12C[M_j]:M_j-\omega(\sqrt{h_j})|M_j|^2.
\]
On its complement we use $W\ge0$. Since the quadratic form
vanishes at zero,
\[
 \frac1{h_j^2}\int_B W(I+h_jM_j)\dd
 \ge\frac12\int_B C[N_j]:N_j\dd-b\omega(\sqrt{h_j}).
\]
The last error tends to zero. Convexity of the quadratic form,
proved in Lemma \ref{lem:quadratic}, gives weak lower
semicontinuity of its integral and hence \eqref{eq:matrix-liminf}.
\end{proof}

\begin{proof}[Proof of Theorem \ref{thm:main}(2)]
Suppose $v_j\to v$ in $L^2$. If the liminf of the energies is
$+\infty$, there is nothing to prove. Otherwise choose a subsequence
attaining that liminf along which the energies are bounded above.
The coercive estimate implies $v_j\in\Hz$ and bounded. By \eqref{eq:compactness}, there exists $v\in\Hz$ and a subsequence (not relabeled) such that 
$D_\rho v_j\rightharpoonup D_\rho v$ in $L^2$.

 Now, having in mind that $\rho$ is compactly supported and therefore $\supp{D_\rho v}\subset \Omega+\supp{\rho}$ for any $v\in H_0^\rho(\Omega;\R^n)$, we apply Lemma \ref{lem:liminf} on a large enough ball
with $M_j=D_\rho v_j$, and use
$C[D_\rho v]:D_\rho v=C[D_{\rho,\mathrm{sym}}v]:D_{\rho,\mathrm{sym}}v$.
The load term converges by strong $L^2$ convergence.
This proves the lower bound.
\end{proof}

\subsection{Recovery sequences}\label{sec:recovery}

\begin{proof}[Proof of Theorem \ref{thm:main}(3)]
Fix $v\in \Hz$ and a sequence $h_j\downarrow0$. By the density result~\eqref{eq:density}, there exists a sequence 
\[ w_k\in C_c^\infty(\Omega;\mathbb{R}^n) \] 
such that 
\[ w_k\to v \quad\text{in }H^\rho(\mathbb{R}^n;\mathbb{R}^n). \] 
For every fixed $k\in\mathbb{N}$, we have 
\[ E_{h_j}(w_k)\to E(w_k) \qquad\text{as }j\to\infty. \] 
Hence, there exists $N_k\in\mathbb{N}$ such that 
\[|E_{h_j}(w_k)-E(w_k)|\leq\frac1k, \text{   for all } j\geq N_k. \] 
We now choose a strictly increasing sequence of indices $(j_k)_k$ recursively by setting 
\[ j_1:=N_1, \qquad j_k:=\max\{N_k,j_{k-1}+1\} \quad\text{for }k\geq2. \] 
It then follows that 
\[ |E_{h_j}(w_k)-E(w_k)|\leq\frac1k, \text{   for all } j\geq j_k . \] 
We define the recovery sequence blockwise by 
\[ v_j:=w_k \qquad\text{whenever }j_k\leq j<j_{k+1}. \] 
The finitely many terms with $j<j_1$ may be chosen arbitrarily in $C_c^\infty(\Omega;\mathbb{R}^n)$. Let $k(j)$ denote the unique integer such that 
\[ j_{k(j)}\leq j<j_{k(j)+1}. \] 
Then $v_j=w_{k(j)}$. Moreover, $k(j)\to\infty$ as $j\to\infty$. Indeed, for every $K\in\mathbb{N}$, 
\[ j\geq j_K \quad\Longrightarrow\quad k(j)\geq K. \] 
Since $w_k\to v$ in $H_\rho(\mathbb{R}^n;\mathbb{R}^n)$, we conclude that 
\[ v_j=w_{k(j)}\to v \quad\text{in }H_\rho(\mathbb{R}^n;\mathbb{R}^n). \] 
Finally, since $j\geq j_{k(j)}$, we have 
\[ \begin{aligned} |E_{h_j}(v_j)-E(v)| &\leq |E_{h_j}(w_{k(j)})-E(w_{k(j)})| + |E(w_{k(j)})-E(v)| \\ &\leq \frac1{k(j)} + |E(w_{k(j)})-E(v)|. \end{aligned} \] 
Both terms on the right-hand side converge to zero because $k(j)\to\infty$ and $E(w_k)\to E(v)$. Therefore, 
\[ E_{h_j}(v_j)\to E(v), \] 
and $(v_j)_j$ is the required recovery sequence.

\end{proof}

\begin{Obse}
For a general $v\in \Hz$, the constant sequence $v_j=v$ need not be a recovery sequence. Indeed, $D_\rho v$ need not be bounded, while $W$ may grow faster than quadratically or even take the value $+\infty$ away from a neighborhood of $I$. The diagonal construction above overcomes this difficulty and preserves the exterior condition on the displacement.
\end{Obse}

\section{Existence, uniqueness, and convergence of almost minimizers} \label{sec:minimizers}

\begin{Teor}[Convergence of minima and almost minimizers] \label{thm:convergence-minimizers} Assume the hypotheses of Theorem~\ref{thm:main}. Then the limit functional $E$ admits a unique minimizer 
\[ v\in \Hz. \] 
Equivalently, $v$ is the unique weak solution of 
\begin{equation}\label{eq:weak} \int_{\mathbb{R}^n} {C}[D_{\rho,\mathrm{sym}}v] :D_{\rho,\mathrm{sym}}w\,dx = \int_\Omega f\cdot w\,dx \qquad \text{for every } w\in \Hz. \end{equation}
Moreover, 
\[ \lim_{h\downarrow0} \inf_{z\in \Hz}E_h(z) = \min_{z\in \Hz}E(z) = E(v). \] 
More precisely, let $h_j\downarrow0$, let $\varepsilon_j\geq0$ satisfy $\varepsilon_j\to0$, and suppose that $v_j\in \Hz$ satisfies 
\[ E_{h_j}(v_j) \leq \inf_{z\in \Hz}E_{h_j}(z) +\varepsilon_j. \] 
Then \[ v_j\to v \quad\text{strongly in } L^2(\mathbb{R}^n;\mathbb{R}^n), \] 
\[ v_j\rightharpoonup v \quad\text{weakly in } \Hz, \] 
and 
\[ E_{h_j}(v_j)\to E(v). \] 
\end{Teor}

\begin{proof}
Functional $E$ admits a unique minimizer since it is quadratic, strictly convex and coercive. Self-adjointness and \eqref{eq:quadratic-bounds} show that
\[
 a(v,w)=\int_{\R^n}C[\Ds v]:\Ds w\dd,
\]
is symmetric and bounded on $\Hz$. Korn and Poincar\'{e} imply
$a(v,v)\ge c\|v\|_{\rho,2}^2$ Therefore, Lax--Milgram's lemma gives a unique solution of
\eqref{eq:weak}, and completing the quadratic form shows that it
is the unique minimizer of $E$.

The uniform equicoercivity bound \eqref{eq:coercivity}, together with
$E_h(0)=0$, gives finite infima of $E_h$ uniformly bounded from above and below.
Recovery for the minimizer of $E$ shows
\[
 \limsup_{h\downarrow0}\inf E_h\le E(v).
\]
Choose almost minimizers along any sequence $h_j\downarrow0$.
Equicoercivity gives a subsequential strong $L^2$ limit $w$,
and the lower bound gives
\[
 E(v)\le E(w)\le\liminf_j E_{h_j}(v_j)=\liminf_j\inf E_{h_j}.
\]
Here a subsequence attaining the liminf is selected before compactness.
This proves convergence of infima and identifies every cluster
point with the unique minimizer $v$. Uniqueness and compactness
give strong $L^2$ convergence of the whole sequence. Boundedness
in $\Hz$ and identification of every weak cluster point give its
weak convergence there. The energy convergence follows as well.
\end{proof}
\section*{Acknowledgements}
We would like to thank C. Kreisbeck, C. Mora-Corral, and F. Seifert for stimulating discussions and for kindly sharing their work with us.
This work was supported by {\it Agencia Estatal de Investigación} (Spain) through grant PID2023-151823NB-I00 and {\it Junta de Comunidades de Castilla-La Mancha} (Spain) through grant SBPLY/23/180225/000023. G.G.-S. is supported by a Doctoral Fellowship by Universidad de Castilla-La Mancha \text{2024-UNIVERS-12844-404}.

\section*{Conflicts of interest}
The authors declare that there are no conflicts of interest regarding the publication of this paper.

\end{document}